\documentclass[10pt]{amsart}

\usepackage{amssymb,amsmath,amsthm,amsfonts}
\usepackage[all]{xy}

\theoremstyle{plain}
\newtheorem{Theorem}{Theorem}
\newtheorem{Proposition}[Theorem]{Proposition}
\newtheorem{Lemma}[Theorem]{Lemma}

\theoremstyle{definition}
\newtheorem{Remark}[Theorem]{Remark}

\begin{document}
\title{An incoherent simple group}
\author{Diego Rattaggi}
\thanks{Supported by the Swiss National Science Foundation, No.\ PP002--68627}
\address{Universit\'e de Gen\`eve,
Section de math\'ematiques,
2--4 rue du Li\`evre, CP 64, 
CH--1211 Gen\`eve 4, Switzerland}
\email{rattaggi@math.unige.ch}
\date{\today}
\begin{abstract}
We give an example of a finitely presented simple group
containing a finitely generated subgroup which is not 
finitely presented.
\end{abstract}
\maketitle

A group is called \emph{coherent}, 
if every finitely generated subgroup is finitely presented.
The class of coherent groups for example includes 
free groups, surface groups, or $3$-manifold groups 
(\cite{Scott}).
On the other hand, the group $F_2 \times F_2$ is incoherent (\cite{Stallings}),
where $F_k$ denotes the free group of rank $k$.
This also follows from the subsequent result of Grunewald:

\begin{Proposition} (Grunewald \cite[Proposition~B]{Grun}) \label{Grunewald}
Let $F$ be a free group of rank $k \geq 2$, generated by $\{s_1, \ldots, s_k \}$. 
Let $r_1, \ldots, r_m$ be words over
$\{s_1, \ldots, s_k \}^{\pm 1}$ and $R$ their normal closure
$\langle \! \langle  r_1, \ldots, r_m \rangle \! \rangle_{F}$ in $F$.
Let $H$ be the group with presentation $\langle s_1, \ldots, s_k \mid r_1, \ldots, r_m \rangle$
and $\phi$ the canonical epimorphism $\phi : F \to H \cong F/R$.
Let $\overline{F}$ be a free group of rank $k$ generated by
$\{t_1, \ldots, t_k \}$
and $\psi$ the isomorphism $F \to \overline{F}$, mapping $s_i$ to $t_i$,
$i = 1, \ldots, k$.
Let $\overline{H} = \langle t_1, \ldots, t_k \mid \psi(r_1), \ldots, \psi(r_m) \rangle$,
and $\tilde{\psi} : H \to \overline{H}$ the isomorphism induced by $\psi$.
Finally, let
$\overline{\phi}$ be the canonical epimorphism 
$\overline{F} \to \overline{H} \cong 
\overline{F} / \psi(R)$ (see the commutative diagram below for a summary of this notation).
Suppose that $H$ is infinite and $R \ne \{1\}$.
Then the group $\{ (s,t) \in F \times \overline{F} : \tilde{\psi}(\phi(s)) = \overline{\phi}(t) \}$
is a subgroup of $F \times \overline{F}$ generated by the $k+m$ elements
$(s_1,t_1), \ldots, (s_k,t_k), (r_1,1), \ldots, (r_m,1)$,
but it is not finitely presented.
\end{Proposition}
\begin{figure}[h]
\[
\xymatrix{
R = \langle \! \langle  r_1, \ldots, r_m \rangle \! \rangle_{F} \ar@{^{(}->}[r] \ar@{->}[d]_-{\psi|_R}^{\cong}
& F = \langle s_1, \ldots, s_k \rangle \ar@{->>}[r]^-{\phi} \ar@{->}[d]_-{\psi}^{\cong}  
& H \cong F/R \ar@{->}[d]_-{\tilde{\psi}}^{\cong} \\
\psi(R) = \langle \! \langle  \psi(r_1), \ldots, \psi(r_m) \rangle \! \rangle_{\overline{F}} 
\ar@{^{(}->}[r] & \overline{F} = \langle t_1, \ldots, t_k \rangle
\ar@{->>}[r]^-{\overline{\phi}} &
\overline{H} \cong \overline{F} / \psi(R)
}
\]
\caption{The setup for Proposition~\ref{Grunewald}}
\end{figure}

Our strategy will be to construct a finitely presented simple group $\Lambda$ containing a subgroup
isomorphic to $F_2 \times F_2$. Proposition~\ref{Grunewald} then allows us to construct
an explicit subgroup of $\Lambda$ generated by three elements, which is not finitely presented.
We first define a group $\Gamma$ which will contain 
$\Lambda$ as a normal subgroup of index $4$.
Let $\Gamma$ be the group with finite presentation 
$\langle a_1, \ldots, a_6, b_1, \ldots, b_5 \mid R_{\Gamma} \rangle$,
where 
\[
R_{\Gamma} := \left\{ \begin{array}{l l l l l}
  a_1 b_1 a_2^{-1} b_2^{-1}, &a_1 b_2 a_1^{-1} b_1^{-1}, &a_1 b_3 a_2^{-1} b_3^{-1}, 
  &a_1 b_4 a_1^{-1} b_4^{-1}, &a_1 b_5 a_2^{-1} b_5, \\
&&&& \\
  a_1 b_5^{-1} a_4 b_5^{-1}, &a_1 b_3^{-1} a_2^{-1} b_2, &a_1 b_1^{-1} a_2^{-1} b_3, 
  &a_2 b_2 a_2^{-1} b_1^{-1}, &a_2 b_4 a_2^{-1} b_4^{-1}, \\
&&&& \\
  a_2 b_5 a_5^{-1} b_5, &a_3 b_1 a_4^{-1} b_2^{-1}, &a_3 b_2 a_3^{-1} b_1^{-1}, 
  &a_3 b_3 a_4^{-1} b_3^{-1}, &a_3 b_4 a_4 b_4, \\
&&&& \\
  a_3 b_5 a_4 b_4^{-1}, &a_3 b_5^{-1} a_6^{-1} b_5^{-1}, &a_3 b_4^{-1} a_4 b_5, 
  &a_3 b_3^{-1} a_4^{-1} b_2, &a_3 b_1^{-1} a_4^{-1} b_3, \\
&&&& \\  
  a_4 b_2 a_4^{-1} b_1^{-1}, &a_5 b_1 a_5^{-1} b_1^{-1}, &a_5 b_2 a_5 b_3^{-1}, 
  &a_5 b_3 a_6^{-1} b_5, &a_5 b_4 a_5^{-1} b_4^{-1}, \\
&&&& \\
  a_5 b_5 a_6^{-1} b_2^{-1}, &a_5 b_2^{-1} a_6 b_3, &a_6 b_1 a_6^{-1} b_3, 
  &a_6 b_2 a_6^{-1} b_4^{-1}, &a_6 b_4 a_6^{-1} b_1
\end{array} \right\}.
\]

We have found $\Gamma$ using programs written in \textsf{GAP}(\cite{GAP}). 
It is constructed such that simultaneously (compare to the proof of Theorem~\ref{Th})
\begin{itemize}
\item $\Gamma < \mathrm{Aut}(\mathcal{T}_{12}) \times \mathrm{Aut}(\mathcal{T}_{10})$,
where $\mathrm{Aut}(\mathcal{T}_{k})$ denotes the group of automorphisms of the 
$k$-regular tree $\mathcal{T}_{k}$,
\item $\Gamma$ (as well as any finite index subgroup of $\Gamma$) does not have any
non-trivial normal subgroup of infinite index by a theorem of Burger-Mozes (\cite{BMII}), 
\item the subgroup $\langle a_1, a_2, a_3, a_4, b_1, b_2, b_3 \rangle_{\Gamma}$ of $\Gamma$
is not residually finite by a theorem of Wise 
(\cite[Main Theorem~II.5.5]{Wise}),
more precisely the element $a_2 a_1^{-1} a_3 a_4^{-1}$ is contained in each finite
index subgroup of $\langle a_1, a_2, a_3, a_4, b_1, b_2, b_3 \rangle_{\Gamma}$ and 
therefore in each finite index subgroup of $\Gamma$,
\item the normal closure 
$\langle \! \langle  a_2 a_1^{-1} a_3 a_4^{-1} \rangle \! \rangle_{\Gamma}$
has finite index in $\Gamma$,
\item the subgroup $\langle a_5, a_1 a_2^{-1}, b_1, b_4 \rangle_{\Gamma}$ 
is isomorphic to $F_2 \times F_2$.
\end{itemize}
The latter statement will follow from Lemma~\ref{F2F2} below,
using a well-known normal form for elements in $\Gamma$:

\begin{Lemma} \label{BrWi}
(Bridson-Wise \cite[Normal Form Lemma~4.3]{BridsonWise})
Any element $\gamma \in \Gamma$ can be written as 
$\gamma = \sigma_a \sigma_b = \sigma_b ' \sigma_a '$,
where $\sigma_a, \sigma_a '$ are freely reduced words in the subgroup $\langle a_1, \ldots, a_6 \rangle_{\Gamma}$
and $\sigma_b, \sigma_b '$ are freely reduced words in $\langle b_1, \ldots, b_5 \rangle_{\Gamma}$.
The words $\sigma_a, \sigma_a ', \sigma_b, \sigma_b '$ are uniquely determined by $\gamma$.
Moreover, $|\sigma_a| = |\sigma_a '|$ and $|\sigma_b| = |\sigma_b '|$,
where $| \cdot |$ is the word length with respect to the standard generators
$\{a_1, \ldots , a_6, b_1, \ldots , b_5 \}^{\pm 1}$ of $\Gamma$.
\end{Lemma}
Note that Lemma~\ref{BrWi} was proved in \cite{BridsonWise} for a certain class of 
fundamental groups of square complexes covered by a product of trees.
The following two lemmas are a direct consequence of the uniqueness statement in Lemma~\ref{BrWi}.

\begin{Lemma} \label{LemmaFree}
The subgroup $\langle a_1, \ldots, a_6 \rangle_{\Gamma}$ is a free group of rank $6$,
and $\langle b_1, \ldots, b_5 \rangle_{\Gamma}$ is a free group of rank $5$.
\end{Lemma}

\begin{Lemma} \label{F2F2}
Let $a, \tilde{a} \in \langle a_1, \ldots, a_6 \rangle_{\Gamma}$
and $b, \tilde{b} \in \langle b_1, \ldots, b_5 \rangle_{\Gamma}$,
such that $ab = ba$, $a\tilde{b} = \tilde{b}a$, $\tilde{a}b = b\tilde{a}$ and
$\tilde{a}\tilde{b} = \tilde{b}\tilde{a}$.
Then the map 
$\langle a, \tilde{a}, b, \tilde{b} \rangle_{\Gamma} \to 
\langle a, \tilde{a} \rangle_{\Gamma} \times \langle b, \tilde{b} \rangle_{\Gamma}$,
given by $a \mapsto (a,1)$, $\tilde{a} \mapsto (\tilde{a},1)$, 
$b \mapsto (1,b)$, $\tilde{b} \mapsto (1, \tilde{b})$
is an isomorphism of groups.
In particular, if moreover
$\langle a, \tilde{a} \rangle_{\Gamma} \cong F_2$
and $\langle b, \tilde{b} \rangle_{\Gamma} \cong F_2$,
then $\langle a, \tilde{a}, b, \tilde{b} \rangle_{\Gamma} \cong F_2 \times F_2$.
\end{Lemma}

We define our main group $\Lambda$ to be the kernel of the surjective homomorphism of groups
\begin{align}
\varphi: \Gamma &\to \mathbb{Z} / 2\mathbb{Z} \times \mathbb{Z} / 2\mathbb{Z} \notag \\
a_1, \ldots , a_6 &\mapsto (1 + 2\mathbb{Z}, 0 + 2\mathbb{Z}), \notag \\
b_1, \ldots , b_5 &\mapsto (0 + 2\mathbb{Z}, 1 + 2\mathbb{Z}). \notag
\end{align}

\begin{Theorem} \label{Th}
The finitely presented group $\Lambda$ is simple and incoherent.
More precisely, the subgroup 
$\langle a_5^2 b_1^2, a_1 a_2^{-1} b_4^2, a_5^2 \rangle_{\Gamma} < \Lambda$
is not finitely presented.
\end{Theorem}

\begin{proof}
The simplicity of $\Lambda$ follows similarly as in 
\cite[Theorem~3.5]{Rat}, but we recall the main steps in the proof.
First note that $\Gamma$ is the fundamental group of a finite square complex
$X$ having a single vertex called $x$, having $6 + 5$ oriented loops 
(identified with $a_1, \ldots, a_6$, $b_1, \ldots, b_5$) and $6 \cdot 5$ squares
(identified with the relators in $R_{\Gamma}$).
Those $30$ squares are carefully chosen such that several conditions simultaneously hold.
For example, the link of $x$ in $X$ is a complete bipartite graph on $12 + 10$
vertices which correspond to $\{a_1, \ldots, a_6 \}^{\pm 1}$ and
$\{b_1, \ldots, b_5 \}^{\pm 1}$.
As a consequence, the universal covering space $\tilde{X}$ is a product of two regular trees
$\mathcal{T}_{12} \times \mathcal{T}_{10}$, and $\Gamma$ is a subgroup
of $\mathrm{Aut}(\mathcal{T}_{12}) \times \mathrm{Aut}(\mathcal{T}_{10})$.
The local actions on $\mathcal{T}_{12}$ of the projection of $\Gamma$ to
the first factor $\mathrm{Aut}(\mathcal{T}_{12})$ are described by
finite permutation groups $P_h^{(k)}(\Gamma) < S_{12 \cdot 11^{k-1}}$,
where $k \in \mathbb{N}$ (see \cite[Chapter~1]{BMII} or \cite{Rat}).
We compute for $k=1$
\begin{align}
P_h^{(1)}(\Gamma) = \langle 
&(9,10)(11,12),
(1,2)(3,4)(5,6,8),
(1,2)(3,4)(5,8,7)(9,10)(11,12), \notag \\
&(3,9)(4,10),
(1,9,3,6,5,2)(4,12,11,8,7,10)
\rangle = A_{12} \notag
\end{align}
and similarly (taking the projection to the second factor of 
$\mathrm{Aut}(\mathcal{T}_{12}) \times \mathrm{Aut}(\mathcal{T}_{10})$, 
thus getting finite permutation groups $P_v^{(k)}(\Gamma) < S_{10 \cdot 9^{k-1}}$)
\begin{align}
P_v^{(1)}(\Gamma) = \langle
&(1,2)(5,6)(8,10,9),
(1,2,3)(5,6)(9,10),
(1,2)(4,5,6,7)(8,10,9), \notag \\
&(1,2,3)(4,5,6,7)(9,10),
(2,5,6,3)(8,9),
(1,7,9,6,5,8)(2,3,10,4)
\rangle = A_{10}. \notag
\end{align}
By the Normal Subgroup Theorem of Burger-Mozes (\cite{BMII}, see also \cite{Rat}),
using the simplicity and high transitivity of $A_{12}$ and $A_{10}$,
the group $\Gamma$ has no non-trivial normal subgroups of infinite index.
This theorem can also be applied to any finite index subgroup of $\Gamma$, 
in particular to $\Lambda < \Gamma$ which is a subgroup of index $4$ by definition.

Next we have to study the \emph{finite} index subgroups of $\Lambda$. We start with the
subgroup $\langle a_1, a_2, a_3, a_4, b_1, b_2, b_3 \rangle_{\Gamma}$ of $\Gamma$. 
It has a presentation
\begin{align}
\langle  a_1, a_2, a_3, a_4, b_1, b_2, b_3 \mid \,
&a_1 b_1 a_2^{-1} b_2^{-1}, \,
a_1 b_2 a_1^{-1} b_1^{-1}, \,
a_1 b_3 a_2^{-1} b_3^{-1}, \,
a_1 b_3^{-1} a_2^{-1} b_2, \notag \\
&a_1 b_1^{-1} a_2^{-1} b_3, \,
a_2 b_2 a_2^{-1} b_1^{-1}, \,
a_3 b_1 a_4^{-1} b_2^{-1}, \,
a_3 b_2 a_3^{-1} b_1^{-1}, \notag \\
&a_3 b_3 a_4^{-1} b_3^{-1}, \,
a_3 b_3^{-1} a_4^{-1} b_2, \,
a_3 b_1^{-1} a_4^{-1} b_3, \,
a_4 b_2 a_4^{-1} b_1^{-1} \notag
\rangle,
\end{align}
since it is the fundamental group $\pi_1(W,x)$ of a finite square complex
$W$ which is embedded in $X$ by construction. 
This implies that $\pi_1(W,x) < \Gamma = \pi_1(X,x)$.
Observe that the $12$ relators in the presentation of $\pi_1(W,x)$
also appear in the presentation of $\Gamma$,
and that the link of the single vertex (again called $x$) in $W$ is a complete bipartite graph
on $8 + 6$ vertices corresponding to
$\{a_1, a_2, a_3, a_4 \}^{\pm 1}$ and $\{b_1, b_2, b_3 \}^{\pm 1}$.
This group $\pi_1(W,x)$ is not residually finite 
and was introduced exactly for this purpose by Wise in \cite{Wise}
where it is called $\pi_1(D)$.
He showed for example that the element 
$a_2 a_1^{-1} a_3 a_4^{-1}$ is contained in each finite
index subgroup of $\pi_1(W,x)$. Consequently, this element is also contained 
in each finite index subgroup of $\Gamma$.
Since $\langle \! \langle  a_2 a_1^{-1} a_3 a_4^{-1} \rangle \! \rangle_{\Gamma}$
has index $4$ in $\Gamma$,
it follows (see \cite{Rat}) that
\[
\Lambda = \langle \! \langle  a_2 a_1^{-1} a_3 a_4^{-1} \rangle \! \rangle_{\Gamma}
= \bigcap
\{ N \lhd \Gamma : N \text{ has finite index in } \Gamma\},
\]
but the latter group is easily seen to have no proper subgroups of finite index,
hence $\Lambda$ is simple.

We show now that $\Lambda$ is incoherent.
First observe that $a_1 a_2^{-1}$ commutes with $b_1$ in $\pi_1(W,x)$ (and therefore in $\Gamma$),
since
\[
a_1 a_2^{-1} b_1 = a_1 b_2 a_2^{-1} = b_1 a_1 a_2^{-1},
\]
using the square relators $a_2 b_2 a_2^{-1} b_1^{-1}$ and
$a_1 b_2 a_1^{-1} b_1^{-1}$ from the presentation of $\pi_1(W,x)$. 
Moreover, we have forced $X$ to contain certain tori, namely 
$a_1 b_4 a_1^{-1} b_4^{-1}$, 
$a_2 b_4 a_2^{-1} b_4^{-1}$,
$a_5 b_1 a_5^{-1} b_1^{-1}$ and
$a_5 b_4 a_5^{-1} b_4^{-1}$.
This implies that 
\[
a_5 b_1 = b_1 a_5, \, 
a_5 b_4 = b_4 a_5
\text{ and } 
a_1 a_2^{-1} b_4 = b_4 a_1 a_2^{-1}
\]
holds in $\Gamma$. By Lemma~\ref{F2F2} and Lemma~\ref{LemmaFree},
\[
\langle a_5, a_1 a_2^{-1}, b_1, b_4 \rangle_{\Gamma} \cong
\langle a_5, a_1 a_2^{-1} \rangle_{\Gamma} \times \langle b_1, b_4 \rangle_{\Gamma}
\cong F_2 \times F_2,
\]
but this group is not contained in $\Lambda$ (recall the definition of $\Lambda$ 
as kernel of $\varphi$).
Therefore, we take the subgroup
\[
\langle a_5^2, a_1 a_2^{-1}, b_1^2, b_4^2 \rangle_{\Gamma} \cong
\langle a_5^2, a_1 a_2^{-1} \rangle_{\Gamma} \times \langle b_1^2, b_4^2 \rangle_{\Gamma}
\cong F_2 \times F_2,
\]
which is obviously a subgroup of $\Lambda$.
To see that $\langle a_5^2 b_1^2, a_1 a_2^{-1} b_4^2, a_5^2 \rangle_{\Gamma}$
is not finitely presented,
we apply Proposition~\ref{Grunewald}
to the following setting:
$k = 2$, $s_1 = a_5^2$, $s_2 = a_1 a_2^{-1}$, $t_1 = b_1^2$, $t_2 = b_4^2$,
and $m = 1$, $r_1 = s_1$,
that is we take
$F = \langle a_5^2, a_1 a_2^{-1} \rangle_{\Gamma} \cong F_2$
and $\overline{F} = \langle b_1^2, b_4^2 \rangle_{\Gamma} \cong F_2$.
It remains to check that $H$ is infinite and $R \ne \{1\}$,
but this is clear since
$H = \langle s_1, s_2 \mid s_1 \rangle \cong \langle a_1 a_2^{-1} \rangle_{\Gamma} \cong \mathbb{Z}$,
and $R = \langle \! \langle  a_5^2 \rangle \! \rangle_{F}$.
\end{proof}

\begin{Remark}
Note that the group $\Lambda$ can be decomposed as amalgamated products
$F_{9} \ast_{F_{97}} F_{9}$ and $F_{11} \ast_{F_{101}} F_{11}$ 
(see \cite[Proposition~1.4]{RatPhD}), 
in particular $\Lambda$ is torsion-free.
\end{Remark}

\begin{Remark}
It is well-known that the word problem is solvable for any
finitely presented simple group.
In fact, by a theorem of Boone-Higman (\cite{BH}), 
a finitely generated group has solvable word problem
if and only if it can be embedded in a simple subgroup
of a finitely presented group.
However, the \emph{generalized} word problem is not solvable for the simple group $\Lambda$, 
since it contains $F_2 \times F_2$ 
(using a result of Miha\u\i lova \cite{Mih}).
Recall that the generalized word problem is solvable for a group $G$ 
if it is decidable for any element $g \in G$ and any finitely generated subgroup
$H < G$ whether or not $g$ lies in $H$.
It is also known that $\Lambda$ has solvable conjugacy problem (being bi-automatic).
\end{Remark}

\begin{Remark}
We mention two other ways to construct finitely presented incoherent simple groups:
One is to directly embed $F_2 \times F_2$ into a virtually simple group by \cite[Theorem~6.5]{BMII}.
The second one is a finitely presented simple group containing
$\mathrm{GL}_4(\mathbb{Z})$ constructed by E.~Scott (\cite{ScottE}).
It is known (\cite{Grun}) that $\mathrm{SL}_4(\mathbb{Z}) < \mathrm{GL}_4(\mathbb{Z})$ is incoherent.
In fact, if $A,B \in \mathrm{SL}_2(\mathbb{Z})$ generate a free group of rank $2$,
and $E$ denotes the identity matrix in $\mathrm{SL}_2(\mathbb{Z})$, then
\[
\left\langle
\left( \begin{array}{c c}
A & 0 \\
0 & E 
\end{array} \right),
\left( \begin{array}{c c}
B & 0 \\
0 & E 
\end{array} \right),
\left( \begin{array}{c c}
E & 0 \\
0 & A 
\end{array} \right),
\left( \begin{array}{c c}
E & 0 \\
0 & B 
\end{array} \right)
\right\rangle_{\mathrm{SL}_4(\mathbb{Z})}  \cong F_2 \times F_2.
\]
Explicitly, one can take
\[
A = \left( \begin{array}{c c}
1 & 2 \\
0 & 1 
\end{array} \right), \;
B = \left( \begin{array}{c c}
1 & 0 \\
2 & 1 
\end{array} \right).
\]
\end{Remark}


\begin{thebibliography}{99}
 \bibitem{BH}
  W.~W.~Boone and G.~Higman,
  `An algebraic characterization of groups with soluble word problem',
  Collection of articles dedicated to the memory of Hanna Neumann, IX.
  J. Austral. Math. Soc. \textbf{18}(1974), 41--53.
 \bibitem{BridsonWise}
  M.~R.~Bridson and D.~T.~Wise, 
  `$VH$ complexes, towers and subgroups of $F\times F$', 
  Math. Proc. Cambridge Philos. Soc. \textbf{126}(1999), no.~3, 481--497.
 \bibitem{BMII}
  M.~Burger and S.~Mozes,
  `Lattices in product of trees',
  Inst. Hautes \'Etudes Sci. Publ. Math. No. 92 (2001), 151--194.
 \bibitem{GAP}
  The GAP Group, GAP -- Groups, Algorithms, and Programming, Version 4.4; 2004. (http://www.gap-system.org)
 \bibitem{Grun}
  F.~J.~Grunewald, 
  `On some groups which cannot be finitely presented', 
  J. London Math. Soc. (2) \textbf{17}(1978), no. 3, 427--436.
 \bibitem{Mih}
  K.~A.~Miha\u\i lova, 
  `The occurrence problem for direct products of groups', (Russian) 
  Dokl. Akad. Nauk SSSR \textbf{119}(1958), 1103--1105.
 \bibitem{RatPhD}
  D.~Rattaggi,
  `Computations in groups acting on a product of trees: 
  normal subgroup structures and quaternion lattices',
  Ph.D.\ thesis, ETH Z\"urich, 2004.
 \bibitem{Rat} 
  D.~Rattaggi,
  `A finitely presented torsion-free simple group', Preprint, 2004,
  available at arXiv:math.GR/0411546.
 \bibitem{ScottE}
  E.~A.~Scott,
  `The embedding of certain linear and abelian groups in finitely presented simple groups',
  J. Algebra \textbf{90}(1984), no. 2, 323--332.
 \bibitem{Scott}
  G.~P.~Scott, 
  `Finitely generated $3$-manifold groups are finitely presented',
  J. London Math. Soc. (2) \textbf{6}(1973), 437--440. 
 \bibitem{Stallings}
  J.~Stallings, 
  `Coherence of $3$-manifold fundamental groups', 
  S\'eminaire Bourbaki, Vol. 1975/76, 28 \`eme ann\'ee, 
  Exp. No. 481, pp. 167--173. Lecture Notes in Math., Vol. 567, 
  Springer, Berlin, 1977. 
 \bibitem{Wise}
  D.~T.~Wise,
  `Non-positively curved squared complexes, aperiodic tilings,
  and non-residually finite groups',
  Ph.D.\ thesis, Princeton University, 1996.
 \end{thebibliography}
\end{document}